\documentclass[11pt]{amsart}

\usepackage{a4wide,enumerate,color,graphicx}
\usepackage{amsmath,bm}
\usepackage{scrextend} 
\allowdisplaybreaks
\usepackage[normalem]{ulem}
\usepackage[dvipsnames]{xcolor}
\usepackage{setspace}
\usepackage{amsbsy}
\usepackage{hyperref}
\usepackage[overload]{empheq}

\newcommand{\R}{{\mathbb R}} 

\newcommand{\del}{\partial}

\newcommand{\dx}{\textnormal{d}x}

\newtheorem{theorem}{Theorem}

\numberwithin{equation}{section}
\hypersetup{
    colorlinks=true,
    linkcolor=blue,
    filecolor=blue,      
    urlcolor=blue,
}

\begin{document}

\title[Multicomponent systems of gases]{Non-isothermal multicomponent flows with mass diffusion and heat conduction}

\author[S. Georgiadis]{Stefanos Georgiadis}
\address{Computer, Electrical and Mathematical Science and Engineering Division,
King Abdullah University of Science and Technology (KAUST),
Thuwal 23955-6900, Saudi Arabia and Institute for Analysis and Scientific Computing, Vienna University of Technology, Wiedner Hauptstra\ss e 8--10, 1040 Wien, Austria}
\email{stefanos.georgiadis@kaust.edu.sa} 

\author[A. J\"ungel]{Ansgar J\"ungel}
\address{Institute for Analysis and Scientific Computing, Vienna University of Technology, Wiedner Hauptstra\ss e 8--10, 1040 Wien, Austria}
\email{juengel@tuwien.ac.at}

\author[A. E. Tzavaras]{Athanasios E. Tzavaras}
\address{Computer, Electrical and Mathematical Science and Engineering Division,
King Abdullah University of Science and Technology (KAUST),
Thuwal 23955-6900, Saudi Arabia}
\email{athanasios.tzavaras@kaust.edu.sa}

\thanks{The first and second authors acknowledge partial support from the Austrian Science Fund (FWF), grants P33010 and F65.
This work has received funding from the European
Research Council (ERC) under the European Union's Horizon 2020 research and innovation programme, ERC Advanced Grant no.~101018153.}

\begin{abstract}
A type-I model of non-isothermal multicomponent systems of gases describing mass diffusive and heat conductive phenomena is presented.  The derivation of the model and a convergence result among thermomechanical theories in the smooth regime are discussed. Furthermore, the global-in-time existence of weak solutions and the weak-strong uniqueness property are established for the corresponding system with zero barycentric velocity.
\end{abstract}

\date{\today}
%35Q35  	PDEs in connection with fluid mechanics
%76M45  	Asymptotic methods, singular perturbations applied to problems in fluid mechanics
%76N10  	Existence, uniqueness, and regularity theory for compressible fluids and gas dynamics
%76R50  	Diffusion
%76T30  	Three or more component flows
%80A17  	Thermodynamics of continua
\maketitle
%\tableofcontents

%--------------------------------------------------------------
%--------------------------------------------------------------
%--------------------------------------------------------------
%--------------------------------------------------------------
%--------------------------------------------------------------

\subsection{Introduction}

Multicomponent systems of gases are systems composed of several constituents. Due to their ubiquity in nature, their dynamics has raised interest in the mathematical literature. The framework of continuum mechanics and differential equations is employed, according to which each component is modeled as a continuous medium, whose behavior is described by balance laws and constitutive relations. The equations of the model read as: 
\begin{align}%[left=\empheqlbrace]
\label{mass}
    \del_t\rho_i +\textnormal{div}(\rho_iv) &=-\textnormal{div}(\rho_iu_i) \\
\label{mom}
    \del_t(\rho v)+\textnormal{div}(\rho v\otimes v) &=\rho b-\nabla p \\
   \del_t \left(\rho e+\frac{1}{2}\rho v^2\right)+\textnormal{div}\left(\bigg(\rho e+\frac{1}{2}\rho v^2\bigg)v\right) &= \textnormal{div}\left(\kappa\nabla\theta-\sum_{i=1}^n(\rho_ie_i+p_i)u_i\right) \label{en} \\
    &\quad - \textnormal{div}(pv)+\rho r+\rho b\cdot v+\sum_{i=1}^n\rho_ib_i\cdot u_i, \nonumber
\end{align} 
where \eqref{mass} corresponds to the mass balance of the $i$-th component, \eqref{mom} to the momentum, \eqref{en} to the energy balance of the mixture and $i=1,\ldots,n$.

The unknowns are the partial mass densities $\rho_i$, the barycentric velocity of the fluid $v$ and the absolute temperature $\theta$. The system is subject to external fields: 
a field of body forces $\rho_ib_i$, with total body force $\rho b:=\sum_{i=1}^n\rho_ib_i$, and another of radiative heat supplies $\rho_ir_i$, with total heat supply $\rho r=\sum_{i=1}^n\rho_ir_i$. Moreover, $\kappa=\kappa(\rho_1,\dots,\rho_n,\theta)$ is the thermal conductivity. The other quantities are determined as functions of the mass densities and the temperature and are the diffusional velocities $u_i$, the total mass density $\rho:=\sum_{i=1}^n\rho_i$, the partial pressures $p_i$, summing to total pressure $p:=\sum_{i=1}^np_i$, and the partial internal energy densities $\rho_ie_i$, which sum up to the internal energy density $\rho e:=\sum_{i=1}^n\rho_ie_i$. 

The diffusional velocities $u_i$ satisfy the Maxwell-Stefan system \cite{Bot}: 
\begin{equation}%\label{linear-system}
\begin{aligned}%[left=\empheqlbrace] 
-\theta\sum_{j\not=i}b_{ij}\rho_i\rho_j(u_i-u_j)=\epsilon d_i \quad
\label{constr}\mbox{under the constraint}\quad
\sum_{i=1}^n\rho_iu_i =0,
\\
\mbox{where} \quad \qquad d_i=\frac{\rho_i}{\rho}(\rho b-\nabla p)+\rho_i\theta\nabla\frac{\mu_i}{\theta}-\theta(\rho_ie_i+p_i)\nabla\frac{1}{\theta}-\rho_ib_i 
\end{aligned}
\end{equation}
are the generalized forces, $\mu_i$ are the chemical potentials and $b_{ij}$ are positive and symmetric coefficients that depend on $\rho_i,\rho_j$ and $\theta$ and model the binary interactions between the $i$-th and the $j$-th components with a strength that is measured by $\epsilon>0$. For an explanation of the origin of the parameter $\epsilon$, we refer to section \ref{section-deriv}.

The remaining thermodynamic quantities are computed from a set of constitutive relations that describe the material response. Throughout this article, a simple mixture of ideal gases is employed, i.e., the thermodynamics of the $i$-th component is described by a Helmholtz free energy density of the form: \begin{equation}
    \label{free-en}\rho_i\psi_i=\theta\frac{\rho_i}{m_i}\left(\ln\frac{\rho_i}{m_i}-1\right)-c_w\rho_i\theta(\ln\theta-1), \tag{IG}
\end{equation} where $m_i$ are the molar masses and $c_w$ the heat capacity, which for simplicity is assumed to be the same for all components. Given $\rho_i\psi_i$, one computes the chemical potentials as $\mu_i=\frac{\del(\rho_i\psi_i)}{\del\rho_i}$, the partial entropy densities as $\rho_i\eta_i=-\frac{\del(\rho_i\psi_i)}{\del\theta}$, the partial internal energy densities by $\rho_ie_i=\rho_i\psi_i+\rho_i\eta_i\theta$ and the partial pressures by the Gibbs-Duhem relation $p_i=-\rho_i\psi_i+\rho_i\mu_i$. Summing up the partial entropies, the partial internal energies and the partial pressures, we obtain the total entropy density, total internal energy density and total pressure, respectively. Under these relations, system \eqref{mass}-\eqref{en} is closed.

Given equations \eqref{mass}-\eqref{en}, an entropy identity can be derived \cite{GT}: \begin{equation*} %\label{entropy-id}
    \del_t(\rho\eta)+\textnormal{div}(\rho\eta v)=\textnormal{div}\left(\frac{\kappa}{\theta}\nabla\theta-\frac{1}{\theta}\sum_{i=1}^n(\rho_i e_i+p_i-\rho_i\mu_i)u_i\right)+\frac{\kappa}{\theta^2}|\nabla\theta|^2-\frac{1}{\theta}\sum_{i=1}^nu_i\cdot d_i
\end{equation*} 
The last two terms capture entropy production, which according to the second law of thermodynamics must be non-negative. The Clausius-Duhem inequality \begin{equation*} %\label{cl-dh}
    \del_t(\rho\eta)+\textnormal{div}(\rho\eta v)\geq\textnormal{div}\left(\frac{\kappa}{\theta}\nabla\theta-\frac{1}{\theta}\sum_{i=1}^n(\rho_i e_i+p_i-\rho_i\mu_i)u_i\right)
\end{equation*} is a manifestation of the previous statement and has a dual role: For smooth solutions it is used to restrict the form of the constitutive relations, while for weak solutions it is regarded as a criterion of thermodynamic admissibility. 

%--------------------------------------------------------------
%--------------------------------------------------------------

\subsection{Derivation of the model} \label{section-deriv}

In a multicomponent system of $n$ species, one may assume that each of the components is described by its own triplet $(\rho_i,v_i,\theta_i)$, i.e., each species has its own mass density, velocity and temperature. Such a model would contain ample information but the number and complexity of the equations makes it difficult to solve and analyze, and it is (at least) challenging to design experiments able to measure all these quantities. For this reason, simplified models are usually investigated, for example those in which each component is characterized by the triplet $(\rho_i,v_i,\theta)$, i.e., each component has its own mass density and velocity, but the model does not distinguish among different temperatures, and the only temperature involved is that one of the mixture, common for all species.  An even further simplified model that does not distinguish among the different velocities and temperatures, with each constituent described by the triplet $(\rho_i,v,\theta)$, where $v$ is the barycentric velocity of the fluid common for all species. The above models are known as type-III, type-II and type-I, respectively. 

The advantage of type-I models is their simplicity; yet the description of diffusive phenomena would be impossible as different velocities are required for the transportation of mass. Thus, one would like to compensate between the simplicity of a type-I model and the information that a type-II model carries. This counterbalance can be reached if one derives a type-I model via a type-II one, which was originally achieved in \cite{BD} in a more general framework including viscous effects and chemical reactions. The starting point is the type-II model
\begin{align}
  & \del_t\rho_i+\textnormal{div}(\rho_iv)=-\textnormal{div}(\rho_iu_i),
  \quad i=1,\ldots,n, \label{mass-II} \\
  & \del_t (\rho_iv_i)+\textnormal{div}(\rho_iv_i\otimes v_i) \label{mom-II} \\
  &\phantom{xx}{}=\rho_ib_i-\rho_i\nabla\mu_i-\frac{1}{\theta}(\rho_ie_i+p_i-\rho_i\mu_i)\nabla\theta-\theta\sum_{j\not=i}b_{ij}\rho_i\rho_j(v_i-v_j), \nonumber \\
  & \del_t  \left(\rho e+\sum_{i=1}^n\frac{1}{2}\rho_iv_i^2\right) +\textnormal{div}\left(\left(\rho e+\sum_{i=1}^n\frac{1}{2}\rho_iv_i^2\right)v\right) 
  \label{en-II} \\
   &\phantom{xx}{}= \textnormal{div}\left(\kappa\nabla\theta-\sum_{i=1}^n(\rho_ie_i+p_i+\frac{1}{2}\rho_iv_i^2)u_i\right) 
        - \textnormal{div}(pv)+\rho b\cdot v+\rho r+\sum_{i=1}^n\rho_ib_i\cdot u_i. \nonumber
\end{align} 

In the context of a type-II model, the barycentric velocity $v$ is defined as $v=\frac{1}{\rho}\sum_{i=1}^n\rho_iv_i$ and the diffusional velocities as $u_i=v_i-v$. The essential difference between \eqref{mass}-\eqref{en} and \eqref{mass-II}-\eqref{en-II} concerns the momentum balances; namely, in the type-I model, only a single momentum balance is available. It serves as an approximation of the $n$ partial momentum balances of the type-II model, in the sense that terms of order $|u_i|^2$ are ignored (cf. \cite{BD}).

The same derivation was obtained in the isothermal case, excluding viscous effects and chemical reactions in \cite{HJTHigh}, using a Chapman-Enskog expansion, where the type-I model was seen as the high-friction limit of the corresponding type-II. To this extent, the last term in \eqref{mom-II}, which corresponds to the friction term due to the interaction between the components, was rescaled by a factor $1/\epsilon$, where $\epsilon>0$ is a relaxation parameter. By letting $\epsilon\to0$, the partial velocities $v_i$ degenerate to a single velocity, which is the barycentric velocity $v$. The emerging type-I model serves as an $\epsilon^2$ approximation of the corresponding type-II model. The additional information takes  the form of a constrained linear system for determining the diffusional velocities $u_i$ and is the isothermal analogue of system \eqref{constr}. The resulting type-I system contains a single velocity $v$, yet the diffusional velocities $u_i$ carry all the information of the mass-diffusive effects. The derivation of \eqref{mass}-\eqref{constr} as a high-friction limit of \eqref{mass-II}-\eqref{en-II} was done for the non-isothermal case in \cite{GT}.

%--------------------------------------------------------------
%--------------------------------------------------------------

\subsection{Dissipative structure}

As was mentioned above, the energy dissipation needs to be non-negative, as it is essential for the model to be compatible with the second law of thermodynamics. The dissipation \begin{equation*} %\label{diss}
    \mathcal{D}=\frac{1}{\theta^2}\kappa|\nabla\theta|^2-\frac{1}{\theta}\sum_{i=1}^nu_i\cdot d_i
\end{equation*} contains two terms: the first term is the dissipation due to heat conduction, while the second one describes dissipation caused by friction among the components. Keeping in mind that this model should encapsulate three different theories (one which describes only mass-diffusive phenomena when the temperature is kept constant; one that describes only thermal effects when the diffusional velocities vanish; and one which combines both phenomena), it is expected that each dissipation term should be non-negative independent of the other. Indeed, due to the non-negativity of the thermal conductivity $\kappa$, as indicated by Fourier's law of heat conduction, the dissipation due to heat conduction is non-negative and one needs to focus only on the second term.

There are two ways to show the non-negativity of the frictional dissipation: the first one consists of substituting $d_i$ by \eqref{constr} and using the symmetry of the coefficients $b_{ij}$ to deduce \begin{equation}\label{mass-diss1}
    -\frac{1}{\theta}\sum_{i=1}^nu_i\cdot d_i=\frac{1}{2}\sum_{i=1}^n\sum_{j\not=i}^nb_{ij}\rho_i\rho_j|u_i-u_j|^2.
\end{equation} 
The second one consists of the inversion of the constrained linear system \eqref{constr} and the subsequent substitution of $u_i$ in the diffusional dissipation. The latter is a delicate process since system \eqref{constr} is singular and thus the existence of a unique solution is not guaranteed. As the generalized forces $d_i$ satisfy 
$\sum_{i=1}^nd_i=0$, the right-hand side of  \eqref{constr} belongs to the range of the matrix on the left-hand side, and thus the system has infinitely many solutions, from which the one satisfying the linear constraint in \eqref{constr} is selected. This procedure was systematically carried out in \cite{HJT} using the Bott-Duffin generalized inverse (also see \cite{BoDr}), and provides an explicit way of solving for the diffusional velocities $u_i$, which in turn allows for estimates for the unknowns to be obtained. 

After the computation of $u_i$, the diffusional dissipation reads: 

\begin{equation}\label{mass-diss2}
    -\frac{1}{\theta}\sum_{i=1}^nu_i\cdot d_i=\sum_{i=1}^n\sum_{j=1}^nA_{ij}\frac{d_i}{\theta\sqrt{\rho_i}}\cdot\frac{d_j}{\theta\sqrt{\rho_j}}
\end{equation} where $A_{ij}$ is the Bott-Duffin inverse of the matrix on the left-hand side of the linear system. The last expression is a quadratic form generated by the matrix $A_{ij}$, which turns out to be positive semi-definite in a particular subspace related to the constraint \eqref{constr}, again verifying that the diffusional dissipation is non-negative.

%--------------------------------------------------------------
%--------------------------------------------------------------

\subsection{Convergence among theories}

Since the energy dissipation is non-negative, the model fits into the general framework of hyperbolic-parabolic systems, as studied in \cite{CT}, for which the existence of a unique local-in-time strong solution has been proved in \cite{GM}. 

As was mentioned above, system \eqref{mass}-\eqref{constr} contains multiple theories encoded in the choice of the parameters $\epsilon$ and $\kappa$. For instance, the choice $\epsilon=0$, $\kappa\not=0$ corresponds to a theory describing heat-conduction but no mass-diffusion. One would like to investigate whether the strong solution of the system with mass-diffusion and heat conduction converges to the strong solution of the system with heat conduction but no mass-diffusion, obtained by setting $\epsilon=0$. The answer is positive and is summarized in the following theorem from \cite{GT}, where $\mathbb{T}^3$ is the three-dimensional torus:

\begin{theorem}\label{thm1}
    Let $\bar{U}^\kappa$ be a strong solution of system \eqref{mass}-\eqref{constr} neglecting mass-diffusive effects (i.e.\ with $\epsilon=0$) defined on a maximal interval of existence $\mathbb{T}^3\times[0,T^*)$ and let $U^{\epsilon,\kappa}$ be a family of strong solutions of \eqref{mass}-\eqref{constr} defined on $\mathbb{T}^3\times[0,T]$, for some $T < T^*$, which emanate from smooth data $\bar{U}_0^\kappa$, $U_0^{\epsilon,\kappa}$, respectively, and satisfy the uniform bounds \begin{equation} \label{bounds}
        0<\delta\leq\rho_j,\bar{\rho}_j\le M, \quad 0<\delta\le\theta,\bar{\theta}\le M 
    \end{equation} for some $\delta,M>0$. Moreover, assume that $0\le\kappa(\rho_1,\dots,\rho_n,\theta)\le M$. Then, $U^{\epsilon,\kappa}\to\bar{U}^\kappa$ in the relative entropy sense, as $\epsilon\to0$.
\end{theorem} 

Similarly, one can simultaneously let $\epsilon,\kappa\to0$, in order to obtain convergence to the adiabatic theory (cf. \cite{GT}):

\begin{theorem}\label{thm2} 
Let $\bar{U}$ be a strong solution of \eqref{mass}-\eqref{constr} with $\epsilon=\kappa=0$, defined on a maximal interval of existence $\mathbb{T}^3\times[0,T^*)$, and let $U^{\epsilon,\kappa}$ be a family of strong solutions of \eqref{mass}-\eqref{constr} defined on $\mathbb{T}^3\times[0,T]$, for $T<T^*$, emanating from smooth data $\bar{U}_0,U_0^{\epsilon,\kappa}$ respectively. Under the hypotheses of Theorem \ref{thm1}, $U^{\epsilon,\kappa}\to\bar{U}$ in the relative entropy sense, as $\epsilon,\kappa\to0$. 
\end{theorem} 

The proof of Theorems \ref{thm1} and \ref{thm2} is based on the relative entropy \cite{CT,GT,HJT},
\[
  \mathcal{H}(U|\bar{U}) = \int_\Omega\bigg[\frac{1}{2}\rho|v-\bar{v}|^2 + \sum_i\frac{1}{m_i}\bigg(\rho_i\log\frac{\rho_i}{\bar{\rho}_i}
  -(\rho_i-\bar{\rho}_i)\bigg) - c_w\rho\bigg(\log\frac{\theta}{\bar{\theta}}
  +(\theta-\bar{\theta})\bigg)\bigg]\dx, 
\]
which can be seen as a measure of the distance between the solutions of the two systems, namely $(U^{\epsilon,\kappa},\bar{U}^\kappa)$ for Theorem \ref{thm1} and $(U^{\epsilon,\kappa},\bar{U})$ for Theorem \ref{thm2}. The relative entropy identity
\[ \begin{split}
    & \del_t\left(\bar{\theta}\mathcal{H}(U|\bar{U})\right)+\textnormal{div}\Big[v\bar{\theta}\mathcal{H}(U|\bar{U})+(p-\bar{p})(v-\bar{v})+\sum_j\rho_ju_j(\mu_j-\bar{\mu}_j) \\
    & - (\theta-\bar{\theta})\left(\frac{1}{\theta}\kappa\nabla\theta-\frac{1}{\bar{\theta}}\bar{\kappa}\nabla\bar{\theta}\right)+(\theta-\bar{\theta})\sum_j\rho_j\eta_ju_j\Big]+\bar{\theta}\kappa\left|\frac{\nabla\theta}{\theta}-\frac{\nabla\bar{\theta}}{\bar{\theta}}\right|^2 \\
    & - \frac{\bar{\theta}}{\theta}\sum_ju_j\cdot d_j=(\del_t\bar{\theta}+\bar{v}\cdot\nabla\bar{\theta})(-\rho\eta)(U|\bar{U})-p(U|\bar{U})\textnormal{div}\bar{v}-(\eta-\bar{\eta})\rho(v-\bar{v})\cdot\nabla\bar{\theta} \\
    & - \sum_j\nabla\bar{\mu}_j\cdot\left(\frac{\rho_j}{\rho}-\frac{\bar{\rho}_j}{\bar{\rho}}\right)\rho(v-\bar{v})-\rho(v-\bar{v})\nabla\bar{v}\cdot(v-\bar{v})-\sum_j\nabla\bar{\mu}_j\cdot\rho_ju_j \\
    & - \left(\frac{\nabla\theta}{\theta}-\frac{\nabla\bar{\theta}}{\bar{\theta}}\right)\cdot\frac{\nabla\bar{\theta}}{\bar{\theta}}(\bar{\theta}\kappa-\theta\bar{\kappa})-\frac{\theta}{\bar{\theta}}\nabla\bar{\theta}\cdot\sum_j\rho_j\eta_ju_j
\end{split} \]
where
 \begin{align*}
  p(U|\bar{U}) &=p-\bar{p}-\sum_j\bar{p}_{\rho_j}(\rho_j-\bar{\rho}_j)-\bar{p}_\theta(\theta-\bar{\theta})
  \\
  (-\rho\eta)(U|\bar{U}) &=-\rho\eta+\bar{\rho}\bar{\eta}+\sum_j(\bar{\rho}\bar{\eta})_{\rho_j}(\rho_j-\bar{\rho}_j)+(\bar{\rho}\bar{\eta})_\theta(\theta-\bar{\theta}) 
  \end{align*}
is then used to obtain a stability estimate for the difference of the two solutions, by controlling the first five terms on the right-hand side by the relative entropy of the two solutions and absorbing the last three terms by the dissipation on the left-hand side, so that the relative entropy of the two solutions is bounded by the relative entropy of the initial data for all times $0<t<T$. Since the two solutions emanate from the same initial data, they coincide for $t>0$; see \cite{GT} for the full proof.

The assumption that the mass densities in \eqref{bounds} are bounded away from zero can be avoided, at the expense of assuming that the free energy densities $\psi_i$ are in $C^3(\overline{\mathcal{U}})$, where $\mathcal{U}$ is a set in the positive cone $(\mathbb{R}^+)^{n+1}$ with $\overline{\mathcal{U}}$ compact, such that: \[ \mathcal{U}=\{ (\rho_1,\dots,\rho_n,\theta)~:~0<\rho_j,\bar{\rho}_j\le M, ~ 0<\delta\le\rho,\bar{\rho}\le M, ~ 0<\delta\le\theta,\bar{\theta}\le M\}. \] However, in the case of the ideal gas \eqref{free-en}, the presence of the logarithm requires the technical hypothesis that mass densities should avoid vacuum (see \cite[Section 5]{GT} for details).

%--------------------------------------------------------------
%--------------------------------------------------------------

\subsection{Mass and thermal diffusion around zero mean flow}

In the case of zero mean flow, i.e. when the barycentric velocity of the mixture $v$ vanishes, the system reads: 
\begin{align}%[left=\empheqlbrace]
    \label{mass-par}\del_t\rho_i+\textnormal{div}J_i & = 0, \\
    \label{momentum-par}\nabla p & = 0, \\
    \label{energy-par}\del_t(\rho e)+\textnormal{div}J_e & = 0,
\end{align} 
where $u_i$ is the unique solution of \eqref{constr} and the fluxes are given by 
\[ 
  J_i=\rho_iu_i \quad\mbox{and} \quad   J_e=-\kappa\nabla\theta+\sum_{i=1}^n(\rho_ie_i+p_i)u_i. 
\] 
Note that the choice $v=0$ does not make the momentum equation disappear completely; in fact it gives the momentum constraint \eqref{momentum-par}, which makes sure that the system remains consistent with the assumption of zero mean flow, since a non-zero pressure gradient would generate motion, which contradicts the choice $v=0$.

The above system falls into the realm of parabolic problems, in the sense that after a change of variables from the set of prime variables $(\rho_1,\dots,\rho_n,\theta)$ to the set of entropy variables $(\mu_1/\theta,\dots,\mu_n/\theta,-1/\theta)$, the matrix of phenomenological coefficients which relates fluxes and entropy variables, namely the matrix $\mathbb{D}$ such that \[ (J_1,\dots,J_n,J_e)^\top=\mathbb{D} \nabla\left(\frac{\mu_i}{\theta},\dots,\frac{\mu_n}{\theta},-\frac{1}{\theta}\right)^\top \] is positive semi-definite (cf. \cite[Section 2]{GJ}).

System \eqref{mass-par}-\eqref{energy-par} is here solved in a bounded domain $\Omega\subset\R^3$ and is completed by the following boundary and initial conditions: \begin{align}\label{BC}
    & J_i\cdot\nu=0,~J_e\cdot\nu=\lambda(\theta-\theta_0)\quad
    \textnormal{on }\partial\Omega,\ t>0, \\ 
    \label{IC}
    & \rho_i(x,0)=\rho_i^0(x),\ \theta(x,0)=\theta^0(x)
    \quad\textnormal{ in } \Omega,
\end{align} 
where $\lambda\geq0$, $\theta_0>0$, and $\nu$ is the exterior unit normal to the boundary $\del\Omega$. The boundary conditions state that mass cannot enter or exit $\Omega$ through the boundary, while heat exchange is allowed, in a manner proportional (by $\lambda$) to the difference of the temperature of the mixture $\theta$ and the background temperature $\theta_0$.

Even though the Maxwell-Stefan system has been studied extensively in the isothermal case, the only known works in the non-isothermal case concerns the local-in-time existence and uniqueness of classical solutions in \cite{HS} and the global-in-time existence of weak solutions in \cite{HJ}. The goal in \cite{GJ} is to obtain global-in-time weak solutions for the above Maxwell-Stefan-Fourier system, that is compatible with thermodynamics and differs from the model presented in \cite{HJ} in several points, as explained in \cite[Section 1]{GJ}.

\begin{theorem} %\label{theorem-existence}
Let $\Omega\subset\mathbb{R}^3$ be a bounded domain with Lipschitz continuous boundary. Assume that the diffusion coefficients are bounded above and continuous in $(\rho_1,\dots,\rho_n,\theta)$ and the thermal conductivity $\kappa$ is continuous in $(\rho_1,\dots,\rho_n,\theta)$ and satisfies the bounds \begin{equation} \label{kappa}
        c_k(1+\theta^2)\leq\kappa(\rho_1,\dots,\rho_n,\theta)\leq C_k(1+\theta^2)
    \end{equation} for some positive constants $c_k,C_k$ and for all $\theta>0$. If the initial data $\rho_i^0\in L^\infty(\Omega)$ are such that the total mass is bounded away from vacuum and infinity  and $\theta^0\in L^\infty(\Omega)$, with $\inf_\Omega\theta^0>0$, then for every $T>0$ there exists a weak solution of \eqref{mass-par}-\eqref{IC} and \eqref{constr}, satisfying $\rho_i>0$ and $\theta>0$ a.e. in $\Omega_T:=\Omega\times(0,T)$ and having the regularity 
    \begin{align*}
        & \rho_i\in L^\infty(\Omega_T)\cap L^2(0,T;H^1(\Omega))\cap H^1(0,T;(H^2(\Omega))^*), \\
        & \theta\in L^2(0,T;H^1(\Omega))\cap W^{1,16/15}(0,T;(W^{2,16}(\Omega))^*),~\log\theta\in L^2(0,T;H^1(\Omega)).
    \end{align*}
\end{theorem}

The proof of the theorem is based on a suitable regularization and uniform estimates from the regularized entropy (or free energy) inequality. More precisely, we discretize the equations in time by the implicit Euler scheme to avoid any issues regarding the time regularity, transform to the so-called entropy variables, defined by the relative chemical potentials, and add an elliptic higher-order regularization. By construction, the entropy variables yield the positivity of the approximate densities and temperature, while the elliptic regularization gives sufficiently regular solutions, making this transformation rigorous. The approximate problem is solved by a fixed-point argument (Leray--Schauder theorem). Some uniform estimates are derived from a discrete, regularized version of the entropy inequality, acquired after choosing the entropy variables as test functions in the weak formulation of the problem and using the technical assumption \eqref{kappa} as well as the positive semi-definiteness of matrix $A_{ij}$ from \eqref{mass-diss2} in the subspace induced by the constraint \eqref{constr}. Then the de-regularization limit can be performed by using an Aubin--Lions compactness argument. For more details, see \cite{GJ}.

The uniqueness of the local-in-time strong solution was established in \cite{Bot} for the isothermal and in \cite{HS} for the non-isothermal case, but for weak solutions the problem remains open. The most general result in the isothermal case is a weak-strong uniqueness property, i.e., whenever there is a strong solution, any weak solution will coincide with the strong one, and can be found in \cite{HJT}. A similar result but for the non-isothermal system \eqref{mass-par}-\eqref{IC} was proved in \cite{GJ}, in the case when no heat exchange is allowed through the boundary, i.e. $\lambda=0$.

\begin{theorem} \label{theorem-uniqueness1}
Let $U=(\rho_1,\dots,\rho_n,\theta)$ be a weak solution of \eqref{mass-par}-\eqref{IC}, with $\lambda=0$ in \eqref{BC}, and let $\bar{U}=(\bar{\rho}_1,\dots,\bar{\rho}_n,\bar{\theta})$ be a strong solution. Assume that there exist $\delta,M>0$ 
such that the weak solution satisfies
\begin{equation}\label{cond.weak} 
  0<\rho_i\leq M, \quad 0<\theta\leq M,
\end{equation}
and the strong solution satisfies 
\begin{equation}\label{cond.strong} 
  0<\delta\leq\bar\rho_i\leq M, \quad 0<\delta\leq\bar\theta\leq M 
\end{equation} 
as well as \[ \nabla\sqrt{\bar{\rho}_i}\in L^\infty_\textnormal{loc}(\Omega\times(0,T)), \quad \nabla\log\bar\theta\in L^\infty_{\textnormal{loc}}
(\Omega\times(0,T)). \]
Moreover, let the thermal conductivity $\kappa$ be Lipschitz continuous as a function of the temperature,  satisfying \eqref{kappa}. Then, if the initial data $U^0=(\rho_1^0,\dots,\rho_n^0,\theta^0)$ and $\bar{U}^0=(\bar\rho_1^0,\dots,\bar\rho_n^0,\bar\theta^0)$ coincide, the two solutions coincide too, i.e. $U(x,t)=\bar{U}(x,t)$ in $\Omega$, for all $0\leq t<T$.
\end{theorem} 

A problematic aspect of Theorem \ref{theorem-uniqueness1} is the assumption that the mass densities of the strong solution need to be bounded away from vacuum. Not only is this a strong mathematical assumption, but also excludes the case of vanishing concentrations that might occur due to the interaction of the components, if one assumes chemical reactions. Keeping this in mind, one can restate the previous theorem, exchanging the assumption on strictly positive mass densities with an assumption on the finiteness of the diffusional velocities, which is a natural hypothesis, since mass is transported at finite speed.

\begin{theorem} \label{theorem-uniqueness2}
Let $U=(\rho_1,\dots,\rho_n,\theta)$ be a weak solution to \eqref{mass-par}-\eqref{IC}, with $\lambda=0$ in \eqref{BC}, and let $\bar{U}=(\bar{\rho}_1,\dots,\bar{\rho}_n,\bar{\theta})$ be a strong solution. Assume that there exist $\delta,M>0$ such that the weak solution satisfies \eqref{cond.weak} and the strong solution satisfies
\begin{equation} \label{cond.strongII}
    0\leq\bar\rho_i\leq M, \quad 0<\delta\leq\bar\theta\leq M
\end{equation} 
as well as \[ \log\bar{\rho}_i\in H^1_{\textnormal{loc}}(\Omega\times(0,T)), \quad \bar{u}_i\in L^\infty_\textnormal{loc}(\Omega\times(0,T)), \quad \nabla\log\bar\theta\in L^\infty_{\textnormal{loc}}(\Omega\times(0,T)). \] Moreover, let the thermal conductivity $\kappa$ be Lipschitz continuous as a function of the temperature, satisfying \eqref{kappa}. Then, if the initial data $U^0=(\rho_1^0,\dots,\rho_n^0,\theta^0)$ and $\bar{U}^0=(\bar\rho_1^0,\dots,\bar\rho_n^0,\bar\theta^0)$ coincide, the two solutions coincide too, i.e. $U(x,t)=\bar{U}(x,t)$ in $\Omega$, for all $0\leq t<T$.
\end{theorem}

The proof of Theorems \ref{theorem-uniqueness1} and \ref{theorem-uniqueness2} is similar with the one of Theorems \ref{thm1} and \ref{thm2}. In this case, the relative entropy reads
 \[ \mathcal{H}(U|\bar{U}) = \int_\Omega\bigg(\sum_i\frac{\rho_i}{m_i}\log\frac{\rho_i}{\bar{\rho}_i}-\sum_i\frac{\rho_i-\bar{\rho}_i}{m_i}-c_w\rho\log\frac{\theta}{\bar{\theta}}+c_w\rho(\theta-\bar{\theta})\bigg)\dx, 
\] since the barycentric velocity is assumed to be zero and the relative entropy identity is used to obtain a stability estimate for the difference of a weak and a strong solution of the same system, namely \eqref{mass-par}-\eqref{IC}; see \cite{GJ}. The difference between the two versions of the theorem lies in the interpretation of the diffusional dissipation: Theorem \ref{theorem-uniqueness1} requires the inversion of system \eqref{constr} and the subsequent elimination of $u_i$ from the diffusional dissipation, resulting in \eqref{mass-diss2} and requiring more assumptions on the mass densities, while in Theorem \ref{theorem-uniqueness2} one eliminates the generalized forces $d_i$ by using \eqref{constr} and the diffusional dissipation takes the form \eqref{mass-diss1}.

%--------------------------------------------------------------
%--------------------------------------------------------------
%--------------------------------------------------------------
%--------------------------------------------------------------
%--------------------------------------------------------------

\end{document}